\documentclass[10pt,twoside]{article}

\def\firstpage{1}

\usepackage{amsmath}
\usepackage{amsthm}
\usepackage{amssymb}
\usepackage{graphicx}
\usepackage{cite}
\usepackage{epsfig}
\usepackage[normal]{subfigure}
\usepackage{paralist}
\usepackage{caption2}
\usepackage{hyperref}
0.3pt \numberwithin{equation}{section} \catcode`@=11
\def\ps@nk{\def\@oddhead{\vbox{\hbox to \hsize{\footnotesize
\footnotesize{\shorttitle \hfill \thepage}} \vspace{1mm} \hrule
\vspace*{-2mm}}}
\def\@evenhead{\vbox{\hbox to \hsize{\footnotesize \footnotesize{\thepage \hfill
\hfill \shortauthor}} \vspace{1mm} \hrule \vspace*{-2mm}}}
\def\@oddfoot{} \def\@evenfoot{}}
\def\ps@first{\def\@oddhead{\vbox{\vspace*{-5mm}
\hbox to \hsize{\footnotesize{}} \vspace{1mm}
\hbox to \hsize{\footnotesize{
}}\vspace{1mm}\hrule \vspace*{-10mm} \break}}
\def\@evenhead{\vbox{\vspace*{-5mm}
\hbox to \hsize{\footnotesize{}} \vspace{1mm}
\hbox to \hsize{\footnotesize{
}} \vspace{1mm}\hrule \vspace*{-10mm} \break} }
\def\@oddfoot{} \def\@evenfoot{}}

\def\Sec{\@Startsection{section}{1}{\z@}
                                   {-3.5ex \@plus -1ex \@minus -.2ex}%
                                   {2.3ex \@plus.2ex}%
                                   {\normalfont\large\bfseries\boldmath}} % czj
\renewcommand{\@seccntformat}[1]{%
{\csname the#1\endcsname.}\hspace{0.5em}}

\def\@Startsection#1#2#3#4#5#6{%
  \if@noskipsec \leavevmode \fi
  \par
  \@tempskipa #4\relax
  \@afterindenttrue
  \ifdim \@tempskipa <\z@
    \@tempskipa -\@tempskipa \@afterindentfalse
  \fi
  \if@nobreak
    \everypar{}%
  \else
    \addpenalty\@secpenalty\addvspace\@tempskipa
  \fi
  \@ifstar
    {\@ssect{#3}{#4}{#5}{#6}}%
    {\@dblarg{\@Sect{#1}{#2}{#3}{#4}{#5}{#6}}}}
\def\@Sect#1#2#3#4#5#6[#7]#8{%
  \ifnum #2>\c@secnumdepth
    \let\@svsec\@empty
  \else
    \refstepcounter{#1}%
    \protected@edef\@svsec{\large \@seccntformat{#1}\relax}%
  \fi
  \@tempskipa #5\relax
  \ifdim \@tempskipa>\z@
    \begingroup
      #6{%
          \@hangfrom{\hskip #3\relax\@svsec \hskip -1mm}%
          \interlinepenalty \@M #8\@@par}
    \endgroup
    \csname #1mark\endcsname{#7}%
    \addcontentsline{toc}{#1}{%
      \ifnum #2>\c@secnumdepth \else
        \protect\numberline{\csname the#1\endcsname}%
      \fi
      #7}%
  \else
    \def\@svsechd{%
      #6{\hskip #3\relax
      \@svsec #8}%
      \csname #1mark\endcsname{#7}%
      \addcontentsline{toc}{#1}{%
        \ifnum #2>\c@secnumdepth \else
          \protect\numberline{\csname the#1\endcsname}%
        \fi
        #7}}%
  \fi
  \@xsect{#5}}

\def\Subsec{\@StartSubsection{subsection}{2}{\z@}%
                                     {-3.25ex\@plus -1ex \@minus -.2ex}%
                                     {1.5ex \@plus .2ex}%
                                     {\normalfont\bfseries\boldmath}}
\def\@StartSubsection#1#2#3#4#5#6{%
  \if@noskipsec \leavevmode \fi
  \par
  \@tempskipa #4\relax
  \@afterindenttrue
  \ifdim \@tempskipa <\z@
    \@tempskipa -\@tempskipa \@afterindentfalse
  \fi
  \if@nobreak
    \everypar{}%
  \else
    \addpenalty\@secpenalty\addvspace\@tempskipa
  \fi
  \@ifstar
    {\@ssect{#3}{#4}{#5}{#6}}%
    {\@dblarg{\@SubSect{#1}{#2}{#3}{#4}{#5}{#6}}}}
\def\@SubSect#1#2#3#4#5#6[#7]#8{%
  \ifnum #2>\c@secnumdepth
    \let\@svsec\@empty
  \else
    \refstepcounter{#1}%
    \protected@edef\@svsec{\@seccntformat{#1}\relax}%
  \fi
  \@tempskipa #5\relax
  \ifdim \@tempskipa>\z@
    \begingroup
      #6{%
          \@hangfrom{\hskip #3\relax\@svsec\hskip -1.5mm}%
          \interlinepenalty \@M #8\@@par}
    \endgroup
    \csname #1mark\endcsname{#7}%
    \addcontentsline{toc}{#1}{%
      \ifnum #2>\c@secnumdepth \else
        \protect\numberline{\csname the#1\endcsname}%
      \fi
      #7}%
  \else
    \def\@svsechd{%
      #6{\hskip #3\relax
      \@svsec #8}%
      \csname #1mark\endcsname{#7}%
      \addcontentsline{toc}{#1}{%
        \ifnum #2>\c@secnumdepth \else
          \protect\numberline{\csname the#1\endcsname}%
        \fi
        #7}}%
  \fi
  \@xsect{#5}}
\def\list#1#2{\ifnum \@listdepth >5\relax \@toodeep \else \global
\advance \@listdepth\@ne \fi \rightmargin \z@ \listparindent\z@
\itemindent\z@ \csname @list\romannumeral\the\@listdepth\endcsname
\def\@itemlabel{#1}\let\makelabel\@mklab \@nmbrlistfalse #2\relax
\@trivlist \parskip 0pt \parindent\listparindent \advance \linewidth
-\rightmargin \advance\linewidth -\leftmargin
\advance\@totalleftmargin \leftmargin \parshape \@ne
\@totalleftmargin \linewidth \ignorespaces} \catcode`@=12
\newtheoremstyle{theorem}{6pt}{6pt}{\itshape}{}{\bfseries}{.}{.5em}{}
\newtheoremstyle{definition}{6pt}{6pt}{\upshape}{}{\bfseries}{.}{.5em}{}

\allowdisplaybreaks[4]
\theoremstyle{plain}
\newtheorem{thm}{Theorem}[section]

\theoremstyle{definition}

\theoremstyle{remark}
\newtheorem{rem}{Remark}[section]
\DeclareMathOperator{\td}{d}

\def\thebibliography#1{\vspace*{-3mm}
\section*{\leftline{{\Large  {\normalfont {\bfseries References}}}}}
\list{[\arabic{enumi}]} {\settowidth \labelwidth{[#1]} \leftmargin
\labelwidth \advance \leftmargin \labelsep \usecounter{enumi}}
\def\newblock{\hskip .11em plus .33em minus .07em} \sloppy \clubpenalty
4000 \widowpenalty 4000 \sfcode`\.=1000 \relax}

\def\shorttitle{A logarithmically completely monotonic function}%%%%%%%%  Please don't be too long. It is up 40 characters.  %%%%%%%%%%%%%%%%%%
\def\shortauthor{F. Qi,  W.-H. Li}
\title{\Large \bf \boldmath \vspace*{-0cm}
\uppercase{A logarithmically completely monotonic function involving the ratio of gamma functions$^*$} }
\author{\large Feng Qi$^{1,2,3,\dag}$ and Wen-Hui Li$^3$}%%%%%%%%  Here is the full name of every author. %%%%%%%%%%%%%%%%%%
\date{}

\begin{document}
\baselineskip 12pt
\maketitle \vspace{-6.5mm}
\begin{quotation}
{\small {\noindent {\bf Abstract}\ \  In the paper, the authors concisely survey and review some functions involving the gamma function and its various ratios, simply state their logarithmically complete monotonicity and related results, and find necessary and sufficient conditions for a new function involving the ratio of two gamma functions and originating from the coding gain to be logarithmically completely monotonic.

\vskip 2mm \noindent {\bf Keywords}\ \  logarithmically completely monotonic function, ratio, gamma function, coding gain, necessary and sufficient condition.

\vskip 2mm \noindent {\bf MSC(2010)}\ \  Primary 33B15; Secondary 26A48, 65R10.}}%%%%%%%%  Here is the 2010 Mathematics Subject Classification of your paper. %%%%%%%%%%%%%%%%%%
\end{quotation}

\thispagestyle{first}\renewcommand{\thefootnote}{\fnsymbol{footnote}}
\footnotetext{\hspace*{-5mm}
\begin{tabular}{@{}r@{}p{10cm}@{}}
$^\dag$& Corresponding author. Email address: qifeng618@gmail.com (F. Qi)\\
$^1$&Institute of Mathematics, Henan Polytechnic University, Jiaozuo City, Henan Province, 454010, China\\
$^2$&College of Mathematics, Inner Mongolia University for Nationalities, Tongliao City, Inner Mongolia Autonomous Region, 028043, China\\
$^3$&Department of Mathematics, School of Science, Tianjin Polytechnic University, Tianjin City, 300387, China\\
$^*$& The first author was supported in part by the National Natural Science Foundation of China (2014JQ1006).
\end{tabular}}

\vspace{-4mm}

\section{Introduction}

In~\cite[Appendix~B]{IEEE-Lee-Tep-3169}, the function
\begin{equation}\label{coding-gain-h}
h(x)=\frac{(2\sqrt\pi\,)^{1/x}[\Gamma(x+1)]^{1/x}}{[\Gamma(x+1/2)]^{1/x}},
\end{equation}
which originated from a coding gain in~\cite[p.~3171, eq.~(15)]{IEEE-Lee-Tep-3169} and where the classical Euler gamma function $\Gamma$ may be defined for $x>0$ by
\begin{equation*}%\label{egamma}
\Gamma(x)=\int^\infty_0t^{x-1} e^{-t}\td t,
\end{equation*}
was proved to be logarithmically completely monotonic on $(0,\infty)$. An infinitely differentiable function $f$ is said to be logarithmically completely monotonic on an interval $I$ if inequalities
\begin{equation*}
(-1)^k[\ln f(x)]^{(k)}\ge0
\end{equation*}
hold on $I$ for all $k\in\mathbb{N}=\{1,2\dotsc\}$.
For more information about the notion ``the logarithmically completely monotonic function'', please refer to~\cite{Atanassov, CBerg, absolute-mon-simp.tex, compmon2, minus-one, auscm-rgmia, subadditive-qi-guo-jcam.tex, SCM-2012-0142.tex, Open-TJM-2003-Banach.tex} and closely related references therein.
\par
In~\cite{ratio-sqrt-gamma-indonesia.tex, minus-one-JKMS.tex}, the function
\begin{equation}\label{2-gamma-ratio}
\frac{[{\Gamma(x+\alpha+1)}]^{1/(x+\alpha)}}{[{\Gamma(x+1)}]^{1/x}}
\end{equation}
was proved to be logarithmically completely monotonic on $(-1,\infty)$ if and only if $\alpha>0$.
\par
In the preprint~\cite{auscm-rgmia} and its formally published version~\cite{e-gam-rat-comp-mon}, the function
\begin{equation}\label{Qi-Berg-Aust-F}
\frac{[\Gamma(x+1)]^{1/x}}x\biggl(1+\frac1x\biggr)^x
\end{equation}
was verified to be logarithmically completely monotonic on $(0,\infty)$. This result was promptly strengthened in~\cite{CBerg} to the function~\eqref{Qi-Berg-Aust-F} being a Stieltjes transform. A Stieltjes transform is a function $f:(0,\infty)\to[0,\infty)$ which can be written in the form
\begin{equation*}%\label{dfn-stieltjes}
f(x)=\frac{a}x+b+\int_0^\infty\frac1{s+x}{\td\mu(s)},
\end{equation*}
where $a,b$ are nonnegative constants and $\mu$ is a nonnegative measure on $(0,\infty)$ such that
\begin{equation*}
\int_0^\infty\frac1{1+s}\td\mu(s)<\infty.
\end{equation*}
More generally, the inclusions
\begin{equation}\label{S-L-C-relation}
\mathcal{L}[I]\subset\mathcal{C}[I]\quad \text{and}\quad \mathcal{S}\setminus\{0\}\subset\mathcal{L}[(0,\infty)]
\end{equation}
were discovered in~\cite{CBerg, absolute-mon-simp.tex, compmon2, minus-one}, where $\mathcal{S}$, $\mathcal{L}[I]$, and $\mathcal{C}[I]$ denote respectively the set of all Stieltjes transforms, the set of all logarithmically completely monotonic functions on an interval $I$, and the set of all completely monotonic functions on $I$. An infinitely differentiable function $f$ is said to be completely monotonic on an interval $I$ if it satisfies
\begin{equation*}
(-1)^if^{(i)}(x)\ge0
\end{equation*}
on $I$ for all $i\in\{0\}\cup\mathbb{N}$. In the literature, we call~\eqref{S-L-C-relation} Qi-Berg's inclusions.
\par
Some properties of the function $[\Gamma(x+1)]^{1/x}$ and its logarithm can be found in~\cite{minus-one, minus-one-JKMS.tex, Coffey-Debrecen-12} and tightly related references therein.
\par
In~\cite{ATSC-ITSF.tex, note-on-neuman-ITSF-simplified.tex}, the monotonicity of the function
\begin{equation}\label{G(s,t;x)}
G_{s,t}(x)=\frac{[\Gamma(1+tx)]^s}{[\Gamma(1+sx)]^t}
\end{equation}
for $x,s,t\in\mathbb{R}$ such that $1+sx>0$ and $1+tx>0$ with $s\ne t$, and its general form
\begin{equation*}%\label{neuman-funct-gen}
g_{a,b}(x)=\frac{[f(bx)]^a}{[f(ax)]^b}
\end{equation*}
for $ax\in I$ and $bx\in I$, where $a$ and $b$ are two real numbers and $f(x)$ is a positive function on an interval $I$, are investigated. For a much complete survey of this topic, please read~\cite[pp.~73--76, Section~7.6]{bounds-two-gammas.tex}.
\par
Let $s$ and $t$ be two real numbers with $s\ne t$,
$\alpha=\min\{s,t\}$ and $\beta>-\alpha$. For $x\in(-\alpha,\infty)$, define
\begin{equation}\label{h-b-dfn}
h_{\beta}(x)=\begin{cases}
\biggl[\dfrac{\Gamma(\beta+t)}{\Gamma(\beta+s)}\cdot
\dfrac{\Gamma(x+s)}{\Gamma(x+t)}\biggr]^{1/(x-\beta)},&x\ne\beta,\\
\exp[\psi(\beta+s)-\psi(\beta+t)],&x=\beta.
\end{cases}
\end{equation}
For $x\in(0,\infty)$ and $\alpha>0$, let
\begin{equation*}%\label{p-def-sandor-alpha}
p_\alpha(x)=\begin{cases} \biggl[\dfrac{\Gamma(\alpha+1)}{\alpha^\alpha}\cdot
\dfrac{x^x}{\Gamma(x+1)}\biggr]^{1/(\alpha-x)},&x\ne\alpha,\\[1em]
\dfrac{\exp[\psi(\alpha+1)-1]}\alpha,&x=\alpha.
\end{cases}
\end{equation*}
The logarithmically complete monotonicity of the functions $h_{\beta}(x)$ and $p_\alpha(x)$ has been studied in~\cite{sandor-gamma-2-ITSF.tex, sandor-gamma-JKMS.tex}. A special case of the function~\eqref{h-b-dfn} came from problems of traffic flow. For more details, please see~\cite[pp.~63--66, Section~6.5]{bounds-two-gammas.tex}.
\par
For given $y\in(-1,\infty)$ and $\alpha\in(-\infty,\infty)$, let
\begin{equation*}%\label{fun}
h_{\alpha,y}(x)=
\begin{cases}
\dfrac1{(x+y+1)^\alpha}\biggl[\dfrac{\Gamma(x+y+1)}{\Gamma(y+1)}\biggr]^{1/x}, &x\in(-y-1,\infty)\setminus\{0\};\\[0.8em]
\dfrac1{(y+1)^\alpha}\exp[\psi(y+1)],&x=0.
\end{cases}
\end{equation*}
Some inequalities and some necessary and sufficient conditions for the function $h_{\alpha,y}(x)$ and its special cases to be logarithmically completely monotonic were provided in~\cite{Open-TJM-2003-Banach.tex, Open-TJM-2003-Ineq-Ext-JAT.tex, Guo-Qi-TJM-03.tex, Extension-TJM-2003.tex, Ya-Ming-Yu-JMAA-09, Zhao-Chu-Jiang} and many other references listed therein.
\par
Let $s$ and $t$ be two real numbers and $\alpha=\min\{s,t\}$. For
$x\in(-\alpha,\infty)$, define
\begin{equation}\label{psidef}
\Psi_{s,t}(x)=
\begin{cases}
\bigg[\dfrac{\Gamma(x+t)}{\Gamma(x+s)}\bigg]^{1/(t-s)},& s\ne t,\\
e^{\psi(x+s)},&s=t.
\end{cases}
\end{equation}
There have been a large amount of literature devoted to inequalities, logarithmically complete monotonicity, asymptotic expansions, applications, and the like, of functions related to the function~\eqref{psidef}. The work in this field, dating back to 1948, has been surveyed and reviewed in~\cite{bounds-two-gammas.tex, Gautschi-Kershaw-TJANT.tex, Wendel2Elezovic.tex-JIA, Wendel-Gautschi-type-ineq-Banach.tex}. Recently, some new results on inequalities and logarithmically complete monotonicity of functions related to~\eqref{psidef} were obtained in~\cite{SCM-2012-0142.tex, notes-best-simple-open-jkms.tex, BAustMS-5984-RV.tex, Com-Mon-Di-Tri-Div-simp.tex, notes-best-simple-cpaa.tex}.
\par
We may classify newly published papers relating to the gamma and polygamma functions and to the logarithmically complete monotonicity into several groups below.
\begin{enumerate}
\item
Some papers having something to do with the unit ball in $\mathbb{R}^n$ are~\cite{Open-TJM-2003-Banach.tex, notes-best-simple-open-jkms.tex, notes-best-simple-cpaa.tex, unit-ball.tex, mon-funct-gamma-unit-ball.tex, refine-Ivady-gamma-PMD.tex}.
\item
Some papers on asymptotic expansions and complete monotonicity for the gamma function $\Gamma$ or for its ratio such as~\eqref{psidef} are~\cite{Buric-Elezovic-2012-Transform, Convexity2CM.tex, ANLY-D-14-0001.tex, Koumandos-Ruijsenaars, Koumandos-Lamprecht-MC-2013, Koumandos-Pedersen-09-JMAA, vali-Feng-Glob.tex, UPB-1635.tex, Mortici-monoburn.tex, Bukac-Sevli-Gamma.tex, Qi-Springer-2012-Srivastava.tex}.
\item
A series of papers on the complete monotonicity of the function $e^{1/t}-\psi'(t)$ and its variants are~\cite{Yang-Fan-2008-Dec-simp.tex, property-psi-ii-final.tex, Bessel-ineq-Dgree-CM.tex, QiBerg.tex, simp-exp-degree-revised.tex}.
\item
Some papers on the notion ``completely monotonic degree'' and its computation are~\cite{Bessel-ineq-Dgree-CM.tex, simp-exp-degree-revised.tex, psi-proper-fraction-degree-two.tex, Geom-Mean-Deg-One-Guo.tex, PAM-Apr-06-2014-0012.tex}.
\item
Some papers related to the function $[\psi'(x)]^2+\psi''(x)$, its divided difference forms, and their $q$-analogues are \cite{SCM-2012-0142.tex, notes-best-simple-open-jkms.tex, BAustMS-5984-RV.tex, notes-best-simple-cpaa.tex, x-4-di-tri-gamma-p(x)-Slovaca.tex, notes-best-simple-equiv.tex, x-4-q-di-tri-gamma.tex, x-4-di-tri-improve.tex, x-4-di-tri-gamma-upper-lower-combined.tex}.
\item
Some applications of the (logarithmically) complete monotonicity to number theory and mean values are published in~\cite{Geom-Mean-Deg-One-Guo.tex, Norlund-No-CM-JNT.tex, 1st-Sirling-Number-2012.tex, Acta-Sinica-B20120547.tex, MIA-3303.tex, Qi-Zhang-Li.tex}.
\end{enumerate}
\par
After concisely surveying and reviewing, now let us return to the function $h(x)$ in~\eqref{coding-gain-h}. It is easy to see that we may rearrange the function $h(x)$ as the form
\begin{equation*}%\label{coding-gain-h-rearr}
h(x)=\frac{(2\sqrt\pi\,)^{1/x}[\Gamma(x+1)]^{1/x}}{[\Gamma(x+1/2)]^{1/x}}
=\biggl[\frac{2\sqrt\pi\,\Gamma(x+1)}{\Gamma(x+1/2)}\biggr]^{1/x}.
\end{equation*}
It is not difficult to see that the function $h(x)$ is not a special case of any functions in~\eqref{2-gamma-ratio} and~\eqref{G(s,t;x)} to~\eqref{psidef}. On the other hand, the function $h(x)$ may be written in a general form
\begin{equation}\label{habc(x)}
h_{a,b;c}(x)=\biggl[c\frac{\Gamma(x+a)}{\Gamma(x+b)}\biggr]^{1/x},
\end{equation}
where $a,b,c>0$ and $x\in(0,\infty)$.
\par
The aim of this paper is to find necessary and sufficient conditions on $a,b,c$ such that the function $h_{a,b;c}(x)$ or its reciprocal is logarithmically completely monotonic on $(0,\infty)$.
\par
Our main results may be stated as the following theorem.

\begin{thm}\label{h_a,b;c(x)-thm}
Let $a,b,c>0$.
\begin{enumerate}
\item
When $a>b$, the function $h_{a,b;c}(x)$ defined by~\eqref{habc(x)} is logarithmically completely monotonic on $(0,\infty)$ if and only if $c\ge\frac{\Gamma(b)}{\Gamma(a)}$;
\item
When $a<b$, the reciprocal of $h_{a,b;c}(x)$ is logarithmically completely monotonic on $(0,\infty)$ if and only if $c\le\frac{\Gamma(b)}{\Gamma(a)}$.
\end{enumerate}
\end{thm}

\begin{rem}
It is clear that the main result in~\cite[Appendix~B]{IEEE-Lee-Tep-3169} is a special case of Theorem~\ref{h_a,b;c(x)-thm} for $a=1$, $b=\frac12$, and $c=2\sqrt\pi\,>\frac{\Gamma(1/2)}{\Gamma(1)}=\sqrt\pi\,$.
\end{rem}

\section{Proof of Theorem~\ref{h_a,b;c(x)-thm}}

Taking the logarithm of $h_{a,b;c}(x)$ and differentiating give
\begin{equation*}
\ln h_{a,b;c}(x)=\frac{\ln c+\ln\Gamma(x+a)-\ln\Gamma(x+b)}x
\end{equation*}
and
\begin{align*}
[\ln h_{a,b;c}(x)]^{(k)}&=\sum_{i=0}^k\binom{k}{i}\biggl(\frac1x\biggr)^{(k-i)} [\ln c+\ln\Gamma(x+a)-\ln\Gamma(x+b)]^{(i)}\\
&=\frac{(-1)^kk!}{x^{k+1}} [\ln c+\ln\Gamma(x+a)-\ln\Gamma(x+b)]\\
&\quad+\sum_{i=1}^k\binom{k}{i}\frac{(-1)^{k-i}(k-i)!}{x^{k-i+1}} \bigl[\psi^{(i-1)}(x+a)-\psi^{(i-1)}(x+b)\bigr]\\
&=\frac{(-1)^kk!}{x^{k+1}}\Biggl\{\ln c+\ln\Gamma(x+a)-\ln\Gamma(x+b)\\
&\quad+\sum_{i=1}^k\frac{(-1)^ix^i}{i!} \bigl[\psi^{(i-1)}(x+a)-\psi^{(i-1)}(x+b)\bigr]\Biggr\}\\
&\triangleq \frac{(-1)^kk!}{x^{k+1}} H_{a,b;c;k}(x)
\end{align*}
for $k\in\mathbb{N}$. It is easy to see that
$$
H_{a,b;c;k}(0)=\ln c+\ln\Gamma(a)-\ln\Gamma(b)=\ln c-\ln\frac{\Gamma(b)}{\Gamma(a)}
$$
and
\begin{align*}
H_{a,b;c;k}'(x)&=\frac{(-1)^kx^k}{k!} \bigl[\psi^{(k)}(x+a)-\psi^{(k)}(x+b)\bigr]\\
&=\frac{x^k}{k!} \bigl[(-1)^k\psi^{(k)}(x+a)-(-1)^k\psi^{(k)}(x+b)\bigr].
\end{align*}
Since
\begin{equation*}%\label{polygamma}
\psi^{(n)}(z)=(-1)^{n+1}\int_0^{\infty}\frac{t^n}{1-e^{-t}}e^{-zt}\td t
\end{equation*}
for $\Re(z)>0$ and $n\in\mathbb{N}$, see~\cite[p.~260, 6.4.1]{abram}, the function $(-1)^k\psi^{(k)}(x)$ is increasing on $(0,\infty)$ for $k\in\mathbb{N}$. Hence, when $a>b$, the derivative $H_{a,b;c;k}'(x)$ is positive on $(0,\infty)$ for $k\in\mathbb{N}$. This means that, when $a>b$, the function $H_{a,b;c;k}(x)$ is increasing on $(0,\infty)$ for $k\in\mathbb{N}$. Therefore, when $a>b$ and $c\ge\frac{\Gamma(b)}{\Gamma(a)}$, the function $H_{a,b;c;k}(x)$ is positive on $(0,\infty)$ for $k\in\mathbb{N}$. Consequently, when $a>b$ and $c\ge\frac{\Gamma(b)}{\Gamma(a)}$, it follows that
\begin{equation}\label{ineq-k-der}
(-1)^k[\ln h_{a,b;c}(x)]^{(k)}=\frac{k!}{x^{k+1}} H_{a,b;c;k}(x)\ge0
\end{equation}
for $k\in\mathbb{N}$ on $(0,\infty)$, that is, the function $h_{a,b;c}(x)$ defined by~\eqref{habc(x)} is logarithmically completely monotonic.
\par
Similarly, when $a<b$, the derivative $H_{a,b;c;k}'(x)$ is negative and the function $H_{a,b;c;k}(x)$ is decreasing on $(0,\infty)$ for $k\in\mathbb{N}$. Therefore, when $a<b$ and $c\le\frac{\Gamma(b)}{\Gamma(a)}$, the function $H_{a,b;c;k}(x)$ is negative and the inequality~\eqref{ineq-k-der} is reversed on $(0,\infty)$ for $k\in\mathbb{N}$. This implies that the reciprocal of the function $h_{a,b;c}(x)$ defined by~\eqref{habc(x)} is logarithmically completely monotonic.
\par
Conversely, if the function $h_{a,b;c}(x)$ is logarithmically completely monotonic, then, by definition, its first logarithmic derivative is negative, that is,
\begin{equation*}
\ln\biggl[\frac{c \Gamma (a+x)}{\Gamma (b+x)}\biggr]-x [\psi(a+x)-\psi(b+x)]\ge0.
\end{equation*}
Taking $x\to0^+$ in the above inequality yields
\begin{equation*}
\ln\biggl[\frac{c \Gamma (a)}{\Gamma (b)}\biggr]\ge0
\end{equation*}
which implies the required necessary and sufficient conditions.
The proof of Theorem~\ref{h_a,b;c(x)-thm} is complete.

\begin{rem}
It is apparent that the proof of Theorem~\ref{h_a,b;c(x)-thm} is slightly simpler than the proof in~\cite[Appendix~B]{IEEE-Lee-Tep-3169} and that Theorem~\ref{h_a,b;c(x)-thm} generalizes the result obtained in~\cite[Appendix~B]{IEEE-Lee-Tep-3169}.
\end{rem}

\begin{rem}
The techniques in the proof of Theorem~\ref{h_a,b;c(x)-thm} and~\cite[Appendix~B]{IEEE-Lee-Tep-3169} were ever appeared in~\cite{minus-one-JKMS.tex}.
\end{rem}

\begin{rem}
This paper is a slightly modified version of the preprint~\cite{gamma-sqrt-ratio.tex}.
\end{rem}

\end{document}